%% file: ms.tex
\documentclass[11pt]{article}
\usepackage[utf8]{inputenc}
\usepackage{graphicx}

\usepackage{amsmath,amssymb,amsfonts}
\usepackage{caption}
\usepackage{subcaption}
\usepackage{color}

\PassOptionsToPackage{hyphens}{url}
\usepackage{hyperref}

\usepackage{amsthm}

\usepackage{ifthen}
\newboolean{showfixes}
\setboolean{showfixes}{true}
\newcommand{\fixes}[1]{\ifthenelse{\boolean{showfixes}}
	{\textcolor{blue}{#1}}{}}
\newboolean{showold}
\setboolean{showold}{false}
\newcommand{\removed}[1]{\ifthenelse{\boolean{showold}}
	{\textcolor{gray}{#1}}{}}

\usepackage{nomencl}
\makenomenclature
\renewcommand{\nomgroup}[1]{%
\item[\bfseries
\ifthenelse{\equal{#1}{A}}{DERs}{%
\ifthenelse{\equal{#1}{B}}{Grid Services}{%
\ifthenelse{\equal{#1}{C}}{Primal-Dual Control}{%
\ifthenelse{\equal{#1}{D}}{Step-size Adaptation}{}}}}%
]}

\begin{document}

\title{Adaptive Primal-Dual Control for Distributed Energy Resource Management}


\author{Joshua Comden \hspace{0.5in}
Jing Wang \hspace{0.5in}
Andrey Bernstein}
\date{}

\maketitle

\begin{abstract}
    With the increased adoption of distributed energy resources (DERs) in distribution networks, their coordinated control with a DER management system (DERMS) that provides grid services (e.g., voltage regulation, virtual power plant) is becoming more necessary.
    One particular type of DERMS using primal-dual control has recently been found to be very effective at providing multiple grid services among an aggregation of DERs;
    however, the main parameter, the primal-dual step size, must be manually tuned for the DERMS to be effective, which can take a considerable amount of engineering time and labor.
    To this end, we design a simple method that self-tunes the step size(s) and adapts it to changing system conditions.
    Additionally, it gives the DER management operator the ability to prioritize among possibly competing grid services.
    We evaluate the automatic tuning method on a simulation model of a real-world feeder in Colorado with data obtained from an electric utility.
    Through a variety of scenarios, we demonstrate that the DERMS with automatically and adaptively tuned step sizes provides higher-quality grid services than a DERMS with a manually tuned step size.
\end{abstract}



\nomenclature[A, 01]{\(\mathcal{I}\)}{Set of DERs controlled by the DER management operator}
\nomenclature[A, 02]{\(P_{i,t}\)}{Active power injected into the distribution network by DER $i$ at time $t$}
\nomenclature[A, 03]{\(Q_{i,t}\)}{Reactive power injected into the distribution network by DER $i$ at time $t$}
\nomenclature[A, 04]{\(f_{i,t}\)}{Function of cost incurred by DER $i$ at time $t$ for power injection $(P_{i,t},Q_{i,t})$}
\nomenclature[A, 05]{\(\mathcal{Y}_{i,t}\)}{Set of feasible power injections $(P_{i,t},Q_{i,t})$ for DER $i$ at time $t$}

\nomenclature[B, 01]{\(\mathcal{K}\)}{Set of grid services provided by the DER management operator}
\nomenclature[B, 02]{\(\mathcal{J}_k\)}{Set of measurement IDs for grid service $k$}
\nomenclature[B, 03]{\(G_{j,t}\)}{Value of the measurement for $j\in\mathcal{J}_k$ at time $t$}
\nomenclature[B, 04]{\(\underline{G}_{j,t}\)}{Lower bound on the measurement $G_{j,t}$}
\nomenclature[B, 05]{\(\overline{G}_{j,t}\)}{Upper bound on the measurement $G_{j,t}$}
\nomenclature[B, 06]{\(g_{j,t}\)}{Power flow equation that relates the power injections $\{P_{i,t},Q_{i,t}\}_{i\in\mathcal{I}}$ to the measurement $G_{j,t}$}

\nomenclature[C, 01]{\(\underline{D}_{j,t}\)}{Dual variable for the constraint to keep $g_{j,t}$ above $\underline{G}_{j,t}$}
\nomenclature[C, 02]{\(\overline{D}_{j,t}\)}{Dual variable for the constraint to keep $g_{j,t}$ below $\overline{G}_{j,t}$}
\nomenclature[C, 03]{\(H_{i,k,t}^{(P)}\)}{Active power injection direction sent from the the coordinator to DER $i$ with regard to grid service $k$ at time $t$.}
\nomenclature[C, 04]{\(H_{i,k,t}^{(Q)}\)}{Reactive power injection direction sent from the the coordinator to DER $i$ with regard to grid service $k$ at time $t$.}
\nomenclature[C, 05]{\(\alpha\)}{Constant step size for the update of the power injection set points}
\nomenclature[C, 06]{\(\beta\)}{Constant step size for the update of the grid service dual variables}
\nomenclature[C, 07]{\(\nu\)}{Regularization parameter for the power injections}
\nomenclature[C, 08]{\(\epsilon\)}{Regularization parameter for the dual variables}

\nomenclature[D, 01]{\(\alpha_i\)}{Constant step size for the update of the power injection set point of DER $i$}
\nomenclature[D, 02]{\(\alpha_{i,t}\)}{Step size for the update of the power injection set point of DER $i$ associated with time $t$}
\nomenclature[D, 03]{\(\beta_k\)}{Constant step size for the update of the dual variables for grid service $k$}
\nomenclature[D, 04]{\(\beta_{k,t}\)}{Step size for the update of the dual variables for grid service $k$ associated with time $t$}
\nomenclature[D, 05]{\(\widetilde{P}_{i,t+1}\)}{Estimated active power injection set point of DER $i$ at time $t+1$ using step size $\alpha_{i,t}$}
\nomenclature[D, 06]{\(\widetilde{Q}_{i,t+1}\)}{Estimated reactive power injection set point of DER $i$ at time $t+1$ using step size $\alpha_{i,t}$}
\nomenclature[D, 07]{\(\widetilde{\underline{D}}_{j,t+1}\)}{Estimated dual variable for the constraint to keep $g_{j,t}$ above $\underline{G}_{j,t}$ using step size $\beta_{j,t}$}
\nomenclature[D, 08]{\(\widetilde{\overline{D}}_{j,t+1}\)}{Estimated dual variable for the constraint to keep $g_{j,t}$ below $\overline{G}_{j,t}$ using step size $\beta_{j,t}$}
\nomenclature[D, 09]{\(S_{\text{cos}}\)}{Cosine similarity}
\nomenclature[D, 10]{\(S^{(PQ)}_{\text{cos},i,t}\)}{Cosine similarity of the estimated change in the direction of the power injection set point pivoted around time $t$ for DER $i$}
\nomenclature[D, 11]{\(S^{(D)}_{\text{cos},k,t}\)}{Cosine similarity of the estimated change in the direction of the dual variables for grid service $k$}
\nomenclature[D, 12]{\(\underline{S}\)}{Lower threshold for the cosine similarity for which the associated step size is decreased}
\nomenclature[D, 13]{\(\overline{S}\)}{Upper threshold for the cosine similarity for which the associated step size is increased}
\nomenclature[D, 14]{\(\overline{\gamma}\)}{Increase factor for the step sizes}
\nomenclature[D, 15]{\(\underline{\gamma}_i^{(PQ)}\)}{Decrease factor for the step size of DER $i$}
\nomenclature[D, 16]{\(\underline{\gamma}_k^{(D)}\)}{Decrease factor for the step size of grid service $k$}

\printnomenclature

\section{Introduction}
\label{sec:intro}
\input{sections/intro}

\section{Problem Formulation}
\label{sec:problem}
\input{sections/problem}

\section{Adaptive Step Size Tuning}
\label{sec:adaptive_ctr}
\input{sections/adaptive_ctr}

\section{Performance Evaluation}
\label{sec:perf_eval}
\input{sections/perf_eval}

\section{Conclusion}
\input{sections/conclusion}

\section*{Acknowledgments}
This work was authored by the National Renewable Energy Laboratory, operated by Alliance for Sustainable Energy, LLC, for the U.S. Department of Energy (DOE) under Contract No. DE-AC36-08GO28308. Funding provided by U.S. Department of Energy Office of Energy Efficiency and Renewable Energy Solar Energy Technologies Office Award Number TCF-21-25008. The views expressed in the article do not necessarily represent the views of the DOE or the U.S. Government. The U.S. Government retains and the publisher, by accepting the article for publication, acknowledges that the U.S. Government retains a nonexclusive, paid-up, irrevocable, worldwide license to publish or reproduce the published form of this work, or allow others to do so, for U.S. Government purposes.



\input{ms.bbl}
\end{document}

%% file: sections/intro.tex
As distributed energy resources (DERs) become more widely adopted deep within distribution networks, their coordination and control are becoming more important and challenging.
If left uncontrolled, the grid as a whole might not be able to operate reliably as a result of voltage instability~\cite{NERC_ERS_abstract2015,NERC_DER2017,mongavoltage}.
Further, the challenges of DER integration have been pushed into the present in the United States with the issue of Federal Energy Regulatory Commission Order 2222, which requires grid operators to allow DER aggregators to participate in wholesale markets~\cite{cano2020ferc}.
Fortunately, newer DERs have capabilities that can assist in providing grid services, such as frequency support, through careful management and coordination~\cite{NERC_ERS_abstract2015,NERC_DER2017}.

The implementation of a methodology that aggregates, integrates, and manages multiple DERs to provide services to the grid is called a distributed energy resource management system (DERMS).
A DERMS may be operated by a third-party aggregator, the distribution system operator (DSO), or even the transmission system operator (TSO).
Nonetheless, this specification does not affect the framework described in this paper; thus, the entity will simply be referred to as the DER management operator.
It is assumed that the DER management operator has a contract with the DSO, TSO, balancing authority, or market operator to provide grid services by controlling the DERMS~\cite{IEEE_DERMS2021,albertini2022overview,reilly2019integration}.

Many DERMS have been proposed during the past decade to coordinate and control groups of DERs to provide services to the grid.
For example, the United Kingdom is simulating a DERMS for future implementation that aggregates DERs into virtual power plants (VPPs) and uses them to provide aggregated active power set point changes as requested by the higher-level system operator and sets individualized parameters for voltage droop control~\cite{sanz2018enhancing,ahmadi2019uk,davarzani2020coordination}.
See \cite{sarmiento2022applications} for more details on VPPs and physical demonstrations of them around the world.
A small city in France is doing a real-world demonstration of a DERMS to lessen the congestion on its transmission network during critical time intervals by coordinating DER aggregators that control different types of DERs to shift load and generation~\cite{foggia2014nice}.
A real-world demonstration is being done in Korea under different voltage control modes (e.g., constant power factor, constant voltage magnitude) or to keep the voltage magnitude within a fixed range~\cite{park2019introduction}.
On a more conceptual level, the work by \cite{kim2021real} has the DERMS coordinator send power injection reference signals to each DER and then has each DER proportionally deviate from the reference signal based on local voltage and frequency measurements.
A simple DERMS controller was designed by \cite{palahalli2022analysis} using a fuzzy logic centralized controller to determine the power injection set points of the DERs to keep the net power in the network at a set value.

Recently, a new class of DERMS has been developed around the primal-dual control methodology.
It frames the power injection decisions of the DERs as the solution to an optimization problem with the grid services encoded as constraints.
The grid service constraints are then dualized and managed by a centralized DERMS coordinator that sends power injection adjustment control signals to local controllers so that the grid services can be satisfied~\cite{dall2016optimal,wang2021performance,padullaparti2022evaluation,dall2017optimal,wang2018design,gan2021cyber,wang2020performance2,wang2020performance1,wang2021voltage,bernstein2019real}.
In \cite{zhao2022distributed}, the coordinator is further distributed so that each local controller takes partial responsibility in updating a subset of the dual variables corresponding to the grid service measurements it can observe and passing the updates to its neighbors.

A major challenge to implementing any of these DER control schemes in a real power system is that they have many parameters that need to be tuned that depend on the properties of the system (e.g., network topology, decision time granularity), the input behavior (e.g., volatility of solar irradiance), and the priority of the objectives for the operator.
The task of tuning parameters can take a considerable amount of time and effort for engineers to accomplish~\cite{padullaparti2022evaluation,wang2020performance2,wang2020performance1,wang2021voltage,wang2021hardware};
thus, making this step of the implementation process shorter and easier will lower the barrier for effective DERMS deployment.
Beyond tuning, the parameters should self-adapt to changes in the system.
The work by \cite{laaksonen2021towards} argues that making the parameters adaptive is necessary and suggests many high-level options to make various proportional control parameters adaptive to the system.

Because of their significant recent development and excellent demonstrated performance in numerical simulations and hardware-in-the-loop evaluations, primal-dual control DERMS are on their way to commercialization.
Thus, this paper focuses on reducing the cost of its deployment by lessening the amount of time and effort an engineer needs to spend tuning its control parameters. 
The most important parameter is the the step size, which determines how sensitive the primal and dual variables are to measured changes in the system.
The intuition to find a good step size is very rough and inexact.
If the value is too low, then the objectives set by the operator are unable to be met; whereas if the value is too high, it can cause the system to oscillate and become unstable.
Further, the values that are considered too low or too high can vary by magnitudes from one system to another and can change depending on the conditions present. The available theoretical conditions on the step size are typically not useful in practice because they rely on specific system information and are very conservative.

Therefore, this paper aims to design a method to automatically tune the step size to achieve the target control performance and make it adaptive under varying conditions in the grid, such as highly fluctuating solar power generation.
Our design was inspired by a rule-based step size adaptation framework from outside the control theory literature that was originally designed to help speed up the convergence of the primal-dual distributed optimization~\cite{yokota2017efficient}.
The main contributions of the paper can be summarized as:
\begin{enumerate}
    \item We formulate the problem of implementing primal-dual control of DERs as a step size tuning problem of the prescribed grid services and DER power injections. (Section \ref{sec:problem})
    \item We design a simple step size tuning method that self-tunes based on operator-determined priorities and automatically adapts to system conditions.  It is straightforward and easy to implement while significantly reducing the tuning time needed to achieve the targeted performance goals. This paves the path for wider deployment of DERMS in the future.  (Section \ref{sec:adaptive_ctr})
    \item We demonstrate the effectiveness of our method over that of manual tuning through extensive numerical simulations of a real-world feeder under various realistic scenarios. (Section \ref{sec:perf_eval})
\end{enumerate}


%% file: sections/problem.tex
\begin{figure}
    \centering
    \includegraphics[width=0.9\columnwidth]{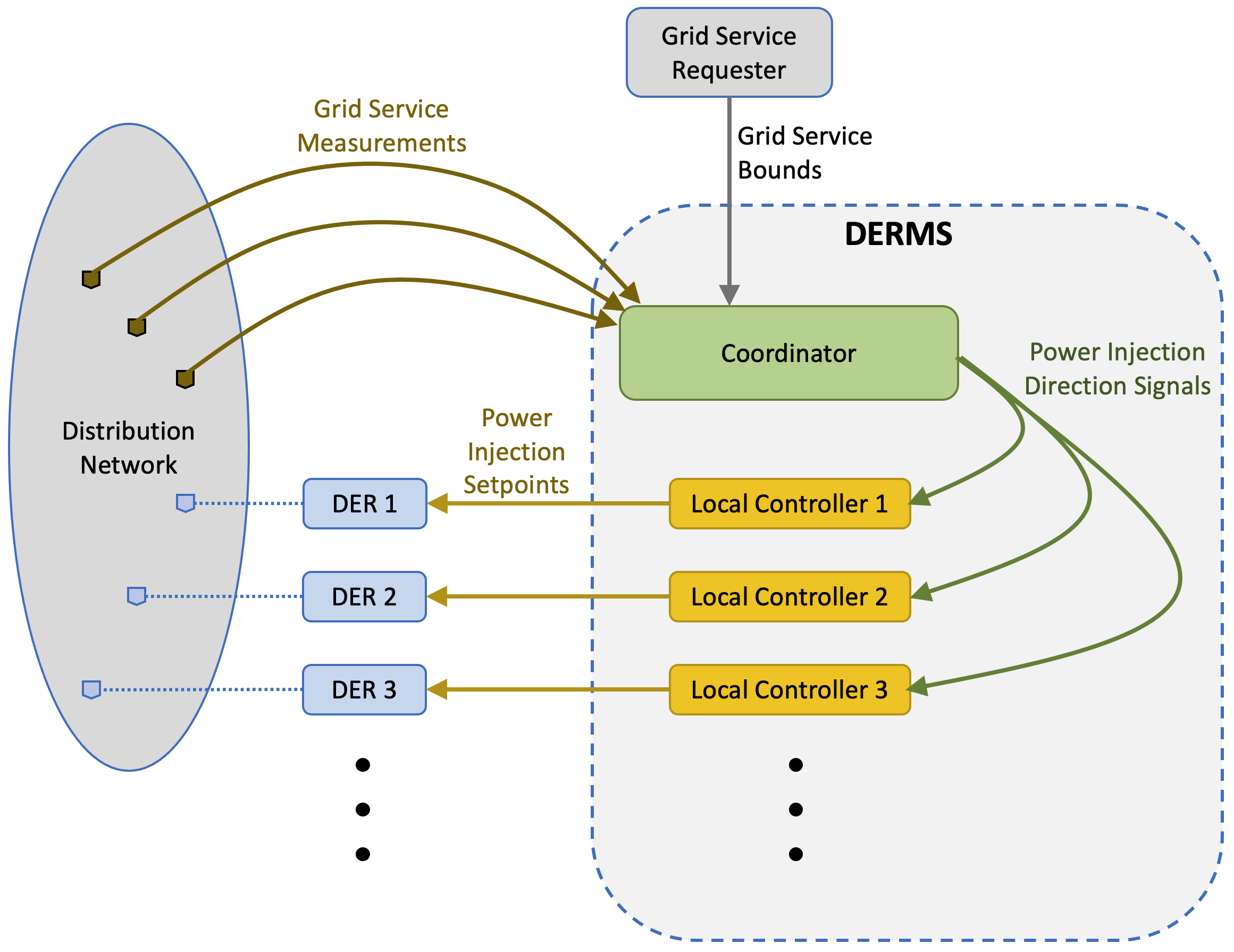}
    \caption{DER management information flow:  The DERMS is tasked with keeping the grid service measurements within the grid service bounds.  The DERMS coordinator collects the grid service information and sends a power injection direction signal to each of its local controllers.  Each local controller uses the signal and local conditions to adjust the power injection set point of its associated DER.}
    \label{fig:DERMS_diagram}
\end{figure}

In this section, we describe a generalized framework of primal-dual control for DER management.
Figure \ref{fig:DERMS_diagram} gives a diagram of the information flow among the system components.
Measurements, communications, and control actions are assumed to be made in equally sized time intervals, $t\in\mathbb{N}$, that can span from subseconds to seconds, depending on the communication infrastructure in place.

\subsection{DER Model}

Let $\mathcal{I}$ be the set of DERs controlled by the DER management operator, where each DER, $i\in\mathcal{I}$, at time $t$, injects $P_{i,t}$, active power, and $Q_{i,t}$, reactive power, into the distribution network.
The power injection $(P_{i,t},Q_{i,t})$ must be within the feasible set, $\mathcal{Y}_{i,t}$, for that DER, which is determined by the type and properties of the DER, such as its state of charge (SOC), inverter rating, and available power.
Detailed concrete examples of feasible sets for various DERs can be found in \cite{bernstein2019real}.

Additionally, the DER incurs a cost of $f_{i,t}(P_{i,t},Q_{i,t})$, which could be the cost of moving away from a desired SOC, the cost of curtailing power from the amount available, the cost of delaying the charging of an electric vehicle, etc.
For example, the cost of the power injection of a photovoltaic (PV) panel can be modeled by the quadratic cost for curtailment, $f_{i,t}(P_{i,t},Q_{i,t})=c^{(P)}_{i}(P_{i,t}-P_{\text{av},i,t})^2 + c^{(Q)}_{i}(Q_{i,t})^2$, where $P_{\text{av},i,t}$ is its maximum available generation, and $(c^{(P)}_{i},c^{(Q)}_{i})$ are weights~\cite{padullaparti2022evaluation}.

\subsection{Grid Service Model}

The DER management operator might have several grid services that it is contracted to provide by controlling the DERs, which could include voltage regulation, limiting power flows over branches, VPPs, etc.~\cite{bernstein2019real}.
Let $\mathcal{K}$ be the set of services provided by the operator.

We define a grid service as a set of measurements from the distribution network and their associated bounds that the operator is tasked with keeping them between.
Specifically, service $k\in\mathcal{K}$ has the set $\mathcal{J}_k$ of measurement IDs, where the value of measurement $j\in\mathcal{J}_k$ at time $t$ is denoted as $G_{j,t}$, with the lower bound, $\underline{G}_{j,t}$, and the upper bound, $\overline{G}_{j,t}$.
For example, if service $k$ is the voltage regulation, then $\{G_{j,t}\}_{j\in\mathcal{J}_k}$ are the voltage magnitude measurements from various locations of the distribution network, and $\{\underline{G}_{j,t},\overline{G}_{j,t}\}_{j\in\mathcal{J}_k}$ are their lower and upper bounds.
If service $k$ is a VPP, then $\{G_{j,t}\}_{j\in\mathcal{J}_k}$ are the power measurements for each phase at the substation, and $\{\underline{G}_{j,t},\overline{G}_{j,t}\}_{j\in\mathcal{J}_k}$ are the time-varying power bands as set by the entity contracted with the DER management operator.

Measurement $G_{j, t}$ is a function of the power injections of the DERs and is described by the power flow equation $g_{j,t}(\{P_{i,t},Q_{i,t}\}_{i\in\mathcal{I}})$, where the uncontrollable power injections in the rest of the distribution network are implicit in the definition of $g_{j,t}$.

\subsection{Primal-Dual Control}

With the feasible sets of the DERs, their costs, and the contracted grid services, the centralized DER management optimization problem to provide the grid services at minimum cost to the DERs is:
\begin{subequations}
    \label{eq:prob_central}
    \begin{align}
        \min_{\{P_{i,t},Q_{i,t}\}_{i\in\mathcal{I}}} \quad & \sum_{i\in\mathcal{I}}f_{i,t}(P_{i,t},Q_{i,t}) \\
        \text{s.t.} \quad & (P_{i,t},Q_{i,t})\in\mathcal{Y}_{i,t}, \quad \forall i\in\mathcal{I} \\
        & \underline{G}_{j,t} \leq g_{j,t}(\{P_{i,t},Q_{i,t}\}_{i\in\mathcal{I}}) \leq \overline{G}_{j,t}, \quad \forall j\in\mathcal{J}_k,\forall k\in\mathcal{K}. \label{eq:prob_central_c}
    \end{align}
\end{subequations}
Let $\{\{\underline{D}_{j,t},\overline{D}_{j,t}\}_{j\in\mathcal{J}_k}\}_{k\in\mathcal{K}}$ be the dual variables for the grid service constraint \eqref{eq:prob_central_c}.

It is impractical, however, for the DER management operator to solve Problem \eqref{eq:prob_central} at every time $t$ and implement the solution for several reasons.
The most difficult and intractable aspect comes from the fact that the power flow equations, $g_{j,t}$, are generally nonconvex and are only partially known because it is assumed that the operator does not have access to measurements of the uncontrollable power injections that are not included in the set of DERs, $\mathcal{I}$.
The other issue is the fact that the operator might not be able to instantaneously receive the current cost function and the feasible set information from every DER under its management.

Instead of solving Problem \eqref{eq:prob_central} directly at each time $t$, primal-dual control employs the use of measurements, dual variables, and gradient approximations to continually move in the direction of the optimal decision trajectory~\cite{bernstein2019online}.
We next describe the primal-dual control DERMS algorithm in a general context. Specific designs for grid services can be found in \cite{dall2016optimal,wang2021performance,padullaparti2022evaluation} for voltage regulation; \cite{dall2017optimal,wang2018design,gan2021cyber,wang2020performance2,wang2020performance1,wang2021voltage} for voltage regulation and VPPs; and \cite{bernstein2019real} for voltage regulation, line current constraints, and VPPs.

The DERMS is structured in a distributed control framework (see Figure \ref{fig:DERMS_diagram}), where a coordinator collects the grid service measurements from the distribution network and sends control signals to a set of local controllers, where each one directly controls a DER or a group of DERs acting as a single DER at a single location in the distribution network.
The control signals can be the common set of dual variables that are sent to all the local controllers~\cite{dall2016optimal,wang2021performance,padullaparti2022evaluation,dall2017optimal,wang2018design,gan2021cyber,bernstein2019real}, or they can be individualized power injection directions~\cite{wang2020performance2,wang2020performance1,wang2021voltage}.
The former option allows the control signal to be a single global message that can be read by each local controller, but the size of the message grows linearly with the number of measurements, and it requires that each local controller has an individualized and up-to-date distribution network model with respect to its power injections.
The latter option, however, requires only the coordinator to hold the whole distribution network model and send individualized small messages to each local controller; this is the option we chose in this paper.
Let $(H_{i,k,t}^{(P)},H_{i,k,t}^{(Q)})$ be the power injection directions sent from the coordinator to DER $i$ with regard to grid service $k$. 

Because the power flow equations, $g_{j,t}$, in general are nonconvex and are not known exactly without error, primal-dual control uses approximated gradients of them with respect to $x\in\{P_{i'},Q_{i'}\}_{i'\in\mathcal{I}}$ evaluated at specific points, denoted as $\widetilde{\nabla}_{(x)} g_{j,t}(\{P_{i',t},Q_{i',t}\}_{i'\in\mathcal{I}})$, to guide the control signals toward the optimal trajectory.
Typically, the approximated gradients are found through linearization of the power flow equations because they remain invariant to changes in the uncontrollable power injections, which allows them to remain unknown for DER management.
As an example, the work by \cite{bernstein2017linear} uses the nodal admittance matrix to find the linearized power flow equations of a multiphase distribution network that relate power injections to voltage magnitudes and power at the substation.

The primal-dual control DERMS algorithm is as follows at each time $t$:
\begin{enumerate}
    \item At the coordinator, for each grid service $k\in\mathcal{K}$:
    \begin{enumerate}
        \item Collect the measurements, $\{G_{j,t}\}_{j\in\mathcal{J}_k}$, from the distribution network and their bounds, $\{(\underline{G}_{j,t},\overline{G}_{j,t})\}_{j\in\mathcal{J}_k}$, from the service requester.
        \item Update the dual variables for each measurement $j\in\mathcal{J}_k$:
        \begin{subequations}
            \label{eq:dual_updates}
            \begin{align}
                \underline{D}_{j,t+1} & := \text{proj}_{\mathbb{R}_+}\left\{\underline{D}_{j,t}+\beta(\underline{G}_{j,t}-G_{j,t}-\epsilon\underline{D}_{j,t})\right\} \\
                \overline{D}_{j,t+1} & := \text{proj}_{\mathbb{R}_+}\left\{\overline{D}_{j,t}+\beta(G_{j,t}-\overline{G}_{j,t}-\epsilon\overline{D}_{j,t})\right\}.
            \end{align}
        \end{subequations}
        \item Send the following power injection direction signals to each $i\in\mathcal{I}$:
        \begin{subequations}
            \begin{align}
                H_{i,k,t+1}^{(P)} & := \sum_{j\in\mathcal{J}_k}(\underline{D}_{j,t+1}-\overline{D}_{j,t+1})\widetilde{\nabla}_{(P_i)} g_{j,t}(\{P_{i',t},Q_{i',t}\}_{i'\in\mathcal{I}}) \\
                H_{i,k,t+1}^{(Q)} & := \sum_{j\in\mathcal{J}_k}(\underline{D}_{j,t+1}-\overline{D}_{j,t+1})\widetilde{\nabla}_{(Q_i)} g_{j,t}(\{P_{i',t},Q_{i',t}\}_{i'\in\mathcal{I}})
            \end{align}
        \end{subequations}
    \end{enumerate}
    \item At each local controller for DER $i\in\mathcal{I}$:
    \begin{enumerate}
        \item Collect the power injection direction signals, $\{(H_{i,k,t+1}^{(P)},H_{i,k,t+1}^{(Q)})\}_{k\in\mathcal{K}}$, from the coordinator.
        \item Update power injection set points for the DER:
        \begin{align}
            \begin{bmatrix}P_{i,t+1}\\Q_{i,t+1}\end{bmatrix} := & \text{proj}_{\mathcal{Y}_{i,t}} \Bigg\{\begin{bmatrix}P_{i,t}\\Q_{i,t}\end{bmatrix} - \alpha\Bigg(\nabla_{(P_i,Q_i)} f_{i,t}(P_{i,t},Q_{i,t}) \nonumber \\
            & \quad\quad\quad\quad\quad\quad -\sum_{k\in\mathcal{K}}\begin{bmatrix}H_{i,k,t+1}^{(P)}\\H_{i,k,t+1}^{(Q)}\end{bmatrix}+\nu\begin{bmatrix}P_{i,t}\\Q_{i,t}\end{bmatrix} \Bigg)\Bigg\} \label{eq:prim_updates}
        \end{align}
    \end{enumerate}
\end{enumerate}
where $(\alpha,\beta,\nu,\epsilon)$ are positive parameters.
Note that the structure of the DERMS coordinator managing the grid services independently allows for multiple coordinators to exist, each taking responsibility for a different grid service.
This could be useful when the measurement data for the different grid services are performed and owned by different entities.

\subsection{Step Size Selection Problem}

The parameters $(\nu,\epsilon)$ and their values are a result of adding regularization terms to the Lagrangian of Problem \eqref{eq:prob_central} in the derivation of the primal-dual control algorithm (very small values are enough; see, e.g., \cite{bernstein2019real}  for details). On the other hand, choosing effective values of the step size parameters $\alpha$ and $\beta$ for the power injection updates \eqref{eq:prim_updates} and dual variable updates \eqref{eq:dual_updates}, respectively, is important but unclear because it depends on the properties and conditions of the distribution network.

If the step size for the power injection updates, $\alpha$, is too small, then the DERs react too slowly to both their own cost and the signals sent from the coordinator, which makes the DERMS ineffective at providing grid services.
On the other hand, if the step size $\alpha$ is too large, then the power injections are too sensitive to the control signals and can overshoot the intended adjustments for grid services.
This can cause oscillations if the coordinator changes the sign of the control signals to counteract the overshoots, which, in turn, causes the power injections to overshoot again but in the opposite direction.
In the best case, the oscillations will dampen out; but in the worst case, they amplify to instability.

The step size for the dual variable updates, $\beta$, similar to the step size for the power injection updates, $\alpha$, is ineffective when it is too small and can cause oscillations when it is too high.
A further complicating fact is that the primal-dual control algorithms in \cite{dall2016optimal,wang2021performance,dall2017optimal,wang2018design,gan2021cyber,wang2020performance2,wang2020performance1,wang2021voltage,bernstein2019real} set the step sizes $\alpha$ and $\beta$ to be equal to each other, whereas only \cite{padullaparti2022evaluation} allows them to be different.
It is possible that the value that is just below the point of causing oscillations for one of them might be too small for the other, which makes the DERMS ineffective at providing grid services.
On the other hand, it is also possible that the value that makes one effective at providing services is too high for the other and cause oscillations.
Thus, the DER management operator will spend a significant amount of time trying to find the balance between being just effective enough at providing grid services while not causing too many oscillations.
The goal of this paper is to enhance the primal-dual control algorithm to automatically find effective and stability-inducing step sizes, which are described in the next section.

%% file: sections/adaptive_ctr.tex
\subsection{Step Size Individualization}

Beyond separating the step sizes $\alpha$ and $\beta$ into the power injection updates \eqref{eq:prim_updates} and the dual variable updates \eqref{eq:dual_updates}, respectively, we further individualize the step size $\alpha$ among the local controllers for each DER to $\{\alpha_i\}_{i\in\mathcal{I}}$ and the step size $\beta$ among the grid services to $\{\beta_k\}_{k\in\mathcal{K}}$.
Even though this greatly expands the number of parameter settings that need to be chosen, we do this for two main reasons.
The first is that it gives the DER management operator the ability to prioritize among grid services, especially when they compete against each other, as will be explained in detail in Section \ref{sec:priority}.
The second is that properties and conditions among the DERs are heterogeneous, so differentiating the step sizes among them allows each local controller to set its step size based on local DER properties and conditions.
Otherwise, if all the local controllers were required to have the same step size, it could burden the communication network to collect all the local condition information at every time $t$ to decide on a common value.

\subsection{Adaptation Procedure}

The adaptation procedure we designed to automatically tune the step sizes is based on  the intuition that values that are too small are ineffective and values that are too large cause oscillations.
The determination is made by comparing the direction that the power injection (dual variables) is heading to the direction it came from, i.e., variations in the time gradients of the decision variables.
If the directions are similar, then it indicates that the power injection (dual variables) has farther to go, and the step size should be increased.
If the directions are very different, however, then it indicates that the power injection (dual variables) is possibly oscillating, and the step size should be decreased.
Let $\{\alpha_{i,t+1}\}_{i\in\mathcal{I}}$ and $\{\beta_{k,t+1}\}_{k\in\mathcal{K}}$ be the time-varying step sizes that are used in setting the power injections $\{P_{i,t+1},Q_{i,t+1}\}_{i\in\mathcal{I}}$ from $\{P_{i,t},Q_{i,t}\}_{i\in\mathcal{I}}$ and the dual variables $\{\{\underline{D}_{j,t+1},\overline{D}_{j,t+1}\}_{j\in\mathcal{J}_k}\}_{k\in\mathcal{K}}$ from $\{\{\underline{D}_{j,t},\overline{D}_{j,t}\}_{j\in\mathcal{J}_k}\}_{k\in\mathcal{K}}$, respectively.

To determine the direction in which the power injections or dual variables are heading before determining the step sizes $\{\alpha_{i,t+1}\}_{i\in\mathcal{I}}$ and $\{\beta_{k,t+1}\}_{k\in\mathcal{K}}$, let $\{\widetilde{P}_{i,t+1},\widetilde{Q}_{i,t+1}\}_{i\in\mathcal{I}}$ and $\{\{\widetilde{\underline{D}}_{j,t+1},\widetilde{\overline{D}}_{j,t+1}\}_{j\in\mathcal{J}_k}\}_{k\in\mathcal{K}}$ be the estimated but not implemented power injection and dual variable settings to act as a placeholder that uses the previous step sizes $\{\alpha_{i,t}\}_{i\in\mathcal{I}}$ and $\{\beta_{k,t}\}_{k\in\mathcal{K}}$, respectively.

We borrow a rule-based step size adaptation framework from \cite{yokota2017efficient}, which uses the cosine similarity and thresholds to determine whether the step sizes should be increased or decreased and by how much.
The cosine similarity, $S_{\text{cos}}(\mathbf{x}_1,\mathbf{x}_2)$, measures how much two same-size vectors $(\mathbf{x}_1,\mathbf{x}_2)$ point in the same or opposite directions with a value between -1 and 1, where -1 indicates complete opposite directions, and 1 indicates the exact same direction:
\begin{align}
    S_{\text{cos}}(\mathbf{x}_1,\mathbf{x}_2) := \frac{\mathbf{x}_1^\top\mathbf{x}_2}{\|\mathbf{x}_1\|_2\|\mathbf{x}_2\|_2}.
\end{align}
The cosine similarity of the estimated change in the direction of the power injection set point pivoted around time $t$ for DER $i\in\mathcal{I}$ is:
\begin{align}
    S^{(PQ)}_{\text{cos},i,t}:=S_{\text{cos}}\left(\begin{bmatrix}\widetilde{P}_{i,t+1}\\\widetilde{Q}_{i,t+1}\end{bmatrix}-\begin{bmatrix}P_{i,t}\\Q_{i,t}\end{bmatrix},\begin{bmatrix}P_{i,t}\\Q_{i,t}\end{bmatrix}-\begin{bmatrix}P_{i,t-1}\\Q_{i,t-1}\end{bmatrix}\right), \label{eq:cosine_sim_prim}
\end{align}
and the cosine similarity of the estimated change in the direction of the dual variables for grid service $k\in\mathcal{K}$ is:
\begin{align}
    S^{(D)}_{\text{cos},k,t}:=S_{\text{cos}}\left(\begin{bmatrix}\widetilde{\underline{\mathbf{D}}}_{k,t+1}\\\widetilde{\overline{\mathbf{D}}}_{k,t+1}\end{bmatrix}-\begin{bmatrix}\underline{\mathbf{D}}_{k,t}\\\overline{\mathbf{D}}_{k,t}\end{bmatrix},\begin{bmatrix}\underline{\mathbf{D}}_{k,t}\\\overline{\mathbf{D}}_{k,t}\end{bmatrix}-\begin{bmatrix}\underline{\mathbf{D}}_{k,t-1}\\\overline{\mathbf{D}}_{k,t-1}\end{bmatrix}\right) \label{eq:cosine_sim_dual}
\end{align}
where $\widetilde{\underline{\mathbf{D}}}_{k,t+1}:=[\underline{D}_{j,t+1}]_{j\in\mathcal{J}_k}$, $\widetilde{\overline{\mathbf{D}}}_{k,t+1}:=[\underline{D}_{j,t+1}]_{j\in\mathcal{J}_k}$, $\underline{\mathbf{D}}_{k,t}:=[\underline{D}_{j,t}]_{j\in\mathcal{J}_k}$, and $\overline{\mathbf{D}}_{k,t}:=[\underline{D}_{j,t}]_{j\in\mathcal{J}_k}$ collects the dual variables into column vectors.
A value of $S^{(PQ)}_{\text{cos},i,t}$ $(S^{(D)}_{\text{cos},k,t})$ close to 1 means that the power injection (dual variables) is heading in a similar direction as the previous time, whereas a value close to -1 means that it is heading in an almost opposite direction.

Finally, the increase or decrease of the step sizes is determined by whether its cosine similarity of the estimated change in direction goes above the threshold, $\overline{S}$, or below the threshold, $\underline{S}$, respectively.
The value of the threshold $\overline{S}$ determines how similar the consecutive directions need to be to indicate a necessary acceleration of that component of the DER management by increasing its step size.
On the other side, the value of the threshold $\underline{S}$ determines how much of a change in direction indicates an oscillation.
The update of the step sizes for the power injections for each DER $i\in\mathcal{I}$ is:
\begin{align}
    \alpha_{i,t+1} := \begin{cases}\overline{\gamma}\alpha_{i,t}, & \text{if}~ S^{(PQ)}_{\text{cos},i,t}> \overline{S} \\ \underline{\gamma}^{(PQ)}_i\alpha_{i,t}, & \text{if}~ S^{(PQ)}_{\text{cos},i,t}< \underline{S} \\
    \alpha_{i,t}, & \text{otherwise}\end{cases} \label{eq:step-size_update_prim}
\end{align}
and for each grid service $k\in\mathcal{K}$ is:
\begin{align}
    \beta_{k,t+1} := \begin{cases}\overline{\gamma}\beta_{k,t}, & \text{if}~ S^{(D)}_{\text{cos},k,t}> \overline{S} \\ \underline{\gamma}^{(D)}_k\beta_{k,t}, & \text{if}~ S^{(D)}_{\text{cos},k,t}< \underline{S} \\
    \beta_{k,t}, & \text{otherwise}\end{cases} \label{eq:step-size_update_dual}
\end{align}
where $\overline{\gamma}>1$ is the common increase factor, and $\{\underline{\gamma}^{(PQ)}_i\}_{i\in\mathcal{I}}$ and $\{\underline{\gamma}^{(D)}_k\}_{k\in\mathcal{K}}$ are the individualized decrease factors that are less than 1.
The principles used to set the values of the factors are further explained in Section \ref{sec:priority}.

Altogether, the primal-dual control DERMS algorithm with the adaptive step size tuning is as follows at each time $t$:
\begin{enumerate}
    \item At the coordinator, for each grid service $k\in\mathcal{K}$:
    \begin{enumerate}
        \item Collect the measurements, $\{G_{j,t}\}_{j\in\mathcal{J}_k}$, from the distribution network and their bounds, $\{(\underline{G}_{j,t},\overline{G}_{j,t})\}_{j\in\mathcal{J}_k}$, from the service requester.
        \item Update the step size $\beta_{k,t+1}$:
        \begin{enumerate}
            \item Estimate the dual variable updates for each measurement, $j\in\mathcal{J}_k$, using the previous step size, $\beta_{k,t}$:
            \begin{subequations}
                \begin{align}
                    \widetilde{\underline{D}}_{j,t+1} & := \text{proj}_{\mathbb{R}_+}\left\{\underline{D}_{j,t}+\beta_{k,t}(\underline{G}_{j,t}-G_{j,t}-\epsilon\underline{D}_{j,t})\right\} \\
                    \widetilde{\overline{D}}_{j,t+1} & := \text{proj}_{\mathbb{R}_+}\left\{\overline{D}_{j,t}+\beta_{k,t}(G_{j,t}-\overline{G}_{j,t}-\epsilon\overline{D}_{j,t})\right\}.
                \end{align}
            \end{subequations}
            \item Calculate the cosine similarity, $S^{(D)}_{\text{cos},k,t}$, of the estimated change in the direction of the dual variables with Equation \eqref{eq:cosine_sim_dual}.
            \item Calculate the step size $\beta_{k,t+1}$ with Equation \eqref{eq:step-size_update_dual}.
        \end{enumerate}
        \item Update the dual variables for each measurement, $j\in\mathcal{J}_k$:
        \begin{subequations}
            \begin{align}
                \underline{D}_{j,t+1} & := \text{proj}_{\mathbb{R}_+}\left\{\underline{D}_{j,t}+\beta_{k,t+1}(\underline{G}_{j,t}-G_{j,t}-\epsilon\underline{D}_{j,t})\right\} \\
                \overline{D}_{j,t+1} & := \text{proj}_{\mathbb{R}_+}\left\{\overline{D}_{j,t}+\beta_{k,t+1}(G_{j,t}-\overline{G}_{j,t}-\epsilon\overline{D}_{j,t})\right\}.
            \end{align}
        \end{subequations}
        \item Send the following power injection direction signals to each $i\in\mathcal{I}$:
        \begin{subequations}
            \begin{align}
                H_{i,k,t+1}^{(P)} & := \sum_{j\in\mathcal{J}_k}(\underline{D}_{j,t+1}-\overline{D}_{j,t+1})\widetilde{\nabla}_{(P_i)} g_{j,t}(\{P_{i',t},Q_{i',t}\}_{i'\in\mathcal{I}}) \\
                H_{i,k,t+1}^{(Q)} & := \sum_{j\in\mathcal{J}_k}(\underline{D}_{j,t+1}-\overline{D}_{j,t+1})\widetilde{\nabla}_{(Q_i)} g_{j,t}(\{P_{i',t},Q_{i',t}\}_{i'\in\mathcal{I}})
            \end{align}
        \end{subequations}
    \end{enumerate}
    \item At each local controller for DER $i\in\mathcal{I}$:
    \begin{enumerate}
        \item Collect the power injection direction signals, $\{(H_{i,k,t+1}^{(P)},H_{i,k,t+1}^{(Q)})\}_{k\in\mathcal{K}}$, from the coordinator.
        \item Update the step size $\alpha_{i,t+1}$:
        \begin{enumerate}
            \item Estimate the power injection set point update using the previous step size, $\alpha_{i,t}$:
            \begin{align}
                \begin{bmatrix}\widetilde{P}_{i,t+1}\\\widetilde{Q}_{i,t+1}\end{bmatrix} := & \text{proj}_{\mathcal{Y}_{i,t}} \Bigg\{\begin{bmatrix}P_{i,t}\\Q_{i,t}\end{bmatrix} - \alpha_{i,t}\Bigg(\nabla_{(P_i,Q_i)} f_{i,t}(P_{i,t},Q_{i,t}) \nonumber \\
                & \quad\quad\quad\quad\quad\quad -\sum_{k\in\mathcal{K}}\begin{bmatrix}H_{i,k,t+1}^{(P)}\\H_{i,k,t+1}^{(Q)}\end{bmatrix}+\nu\begin{bmatrix}P_{i,t}\\Q_{i,t}\end{bmatrix} \Bigg)\Bigg\} 
            \end{align}
            \item Calculate the cosine similarity, $S^{(PQ)}_{\text{cos},i,t}$, of the estimated change in the direction of the power injection with Equation \eqref{eq:cosine_sim_prim}.
            \item Calculate the step size $\alpha_{i,t+1}$ with Equation \eqref{eq:step-size_update_prim}.
        \end{enumerate}
        \item Update the power injection set points for the DER:
        \begin{align}
            \begin{bmatrix}P_{i,t+1}\\Q_{i,t+1}\end{bmatrix} := & \text{proj}_{\mathcal{Y}_{i,t}} \Bigg\{\begin{bmatrix}P_{i,t}\\Q_{i,t}\end{bmatrix} - \alpha_{i,t+1}\Bigg(\nabla_{(P_i,Q_i)} f_{i,t}(P_{i,t},Q_{i,t}) \nonumber \\
            & \quad\quad\quad\quad\quad\quad -\sum_{k\in\mathcal{K}}\begin{bmatrix}H_{i,k,t+1}^{(P)}\\H_{i,k,t+1}^{(Q)}\end{bmatrix}+\nu\begin{bmatrix}P_{i,t}\\Q_{i,t}\end{bmatrix} \Bigg)\Bigg\}.
        \end{align}
    \end{enumerate}
\end{enumerate}

\subsection{DER Management Priority Selection}
\label{sec:priority}

The evolution of the values of the step sizes during the DER management operation are directly affected by the values of the increase factor, $\overline{\gamma}$, and the decrease factors, $\{\underline{\gamma}^{(PQ)}_i\}_{i\in\mathcal{I}}$ and $\{\underline{\gamma}^{(D)}_k\}_{k\in\mathcal{K}}$.
By setting different values among them, we can prioritize the grid services and individually regulate the responsiveness of each DER.

The lesser the decrease factor for a grid service is from 1, the lower priority it has.
This is because it will decrease its step size faster when it experiences oscillations, and thus it become less responsive to its measurements going outside of its bounds compared to other grid services with decrease factors that are closer to 1.
In other words, the adaptation procedure senses possible instability, so it releases some degree of its responsibility to keep its grid services within bounds; the service with the lower decrease factor releases faster.
The DER management operator can set which grid service should be sacrificed first, second, etc., when the coordinator senses instability by the increasing order of the decrease factors.

Similarly, the lesser the decrease factor for a DER is from 1, the faster it will decrease its step size when it experiences local oscillations, and thus it becomes less responsive to the control signals from the coordinator compared to other DERs that have decrease factors that are closer to 1.
This can be used to protect some fragile DERs that can be damaged if their set points are changed too drastically and too frequently by reducing their decrease factors below that of other DERs.

The same increase factor, $\overline{\gamma}$, is used across all the grid services and DERs so that the whole system can accelerate proportionally at the same rate when there are no oscillations present.
This can be set based on the smallest time interval that the DER management operator wants to allow the step sizes to increase by an order of magnitude.

%% file: sections/perf_eval.tex
This section shows through real-world numerical simulations that our step size adaptation procedure for primal-dual control DERMS automatically tunes itself based on system conditions and outperforms manual tuning.

\subsection{Setup}

\begin{figure}
    \centering
    \includegraphics[width=0.7\columnwidth]{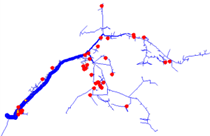}
    \caption{Feeder topology with the participating customers in red.}
    \label{fig:topology}
\end{figure}

Suppose that a DER management operator with a primal-dual control DERMS is tasked with both keeping the voltage magnitudes within prescribed bounds and the active power for each of the three phases at the distribution feeder head between time-varying bounds, i.e., acting as a VPP~\cite{dall2017optimal,wang2018design,gan2021cyber,wang2020performance2,wang2020performance1,wang2021voltage}.
The simulated distribution feeder is based on a real system in Colorado, which contains approximately 2,000 three-phase nodes, with a peak load of 4.6 MW, including 163 customers that also have curtailable PV solar power generation and an inverter; 140 of the 163 customers have energy storage batteries.
Figure \ref{fig:topology} shows the feeder topology with locations of the participating customers.

To evaluate the effectiveness of the adaptive step size tuning procedure, we simulate it under four scenarios that are aimed at investigating different aspects of primal-dual control.
The first scenario looks at the situation of an operator setting the initial step sizes either low or high by a couple of orders of magnitude, which can occur during manual tuning.
The second scenario looks at how changing the VPP set point tracking bounds affects the step sizes, power injection updates, and grid services.
The third scenario looks at how fluctuating PV generation can be dealt with under DER management to satisfy the grid services and how the step sizes adapt to the situation.
The final scenario investigates when the admittance matrix changes through an advanced distribution management system (ADMS) adjusting the load tap changer (LTC) tap position.
This is important because it can create a mismatch between the actual and assumed admittance matrices because changes in the admittance matrix caused by an ADMS are not necessarily immediately communicated to a DER management operator that assumes a specific admittance matrix to do the power flow linearization.

The set of grid services is $\mathcal{K}=(\text{volt},\text{VPP})$ where $\text{volt}$ and $\text{VPP}$ represent voltage regulation and VPP set point tracking, respectively.
Voltage regulation is used to keep the voltage magnitudes of the measured locations between 0.95 p.u. and 1.03 p.u.
The VPP tracks a time-varying set point of the active power injections of each phase at the feeder head with the upper and lower bounds $\pm$10 kW of the set point.

The adaptive step size tuning procedure uses the following parameter settings unless otherwise noted: $\underline{S}=0$, $\overline{S}=0.9$, $\overline{\gamma}=1.005$, $\underline{\gamma}_i=0.95:\forall i\in\mathcal{I}$, $\underline{\gamma}_\text{volt}=0.995$, and $\underline{\gamma}_\text{VPP}=0.5$, which prioritizes voltage regulation over tracking the VPP set point.
The values of the cosine similarity thresholds $(\underline{S},\overline{S})$ are the values suggested by \cite{yokota2017efficient} under a different context.
The primal-dual control with step size tuning is compared against the manually tuned primal-dual control step sizes of $\alpha=\beta=6$, which was found to be the largest value that does not cause undamped oscillations at any local controller.
This value was found by starting from a very low step size and doubling it after each run of the simulation until oscillations were detected.
After that point, the value was decremented by one significant digit until the value of $6$ was found to not cause oscillations.
In fact, this manual tuning process is a very crude version of the automatic tuning process.

\begin{figure}
    \centering
    \includegraphics[width=0.9\columnwidth]{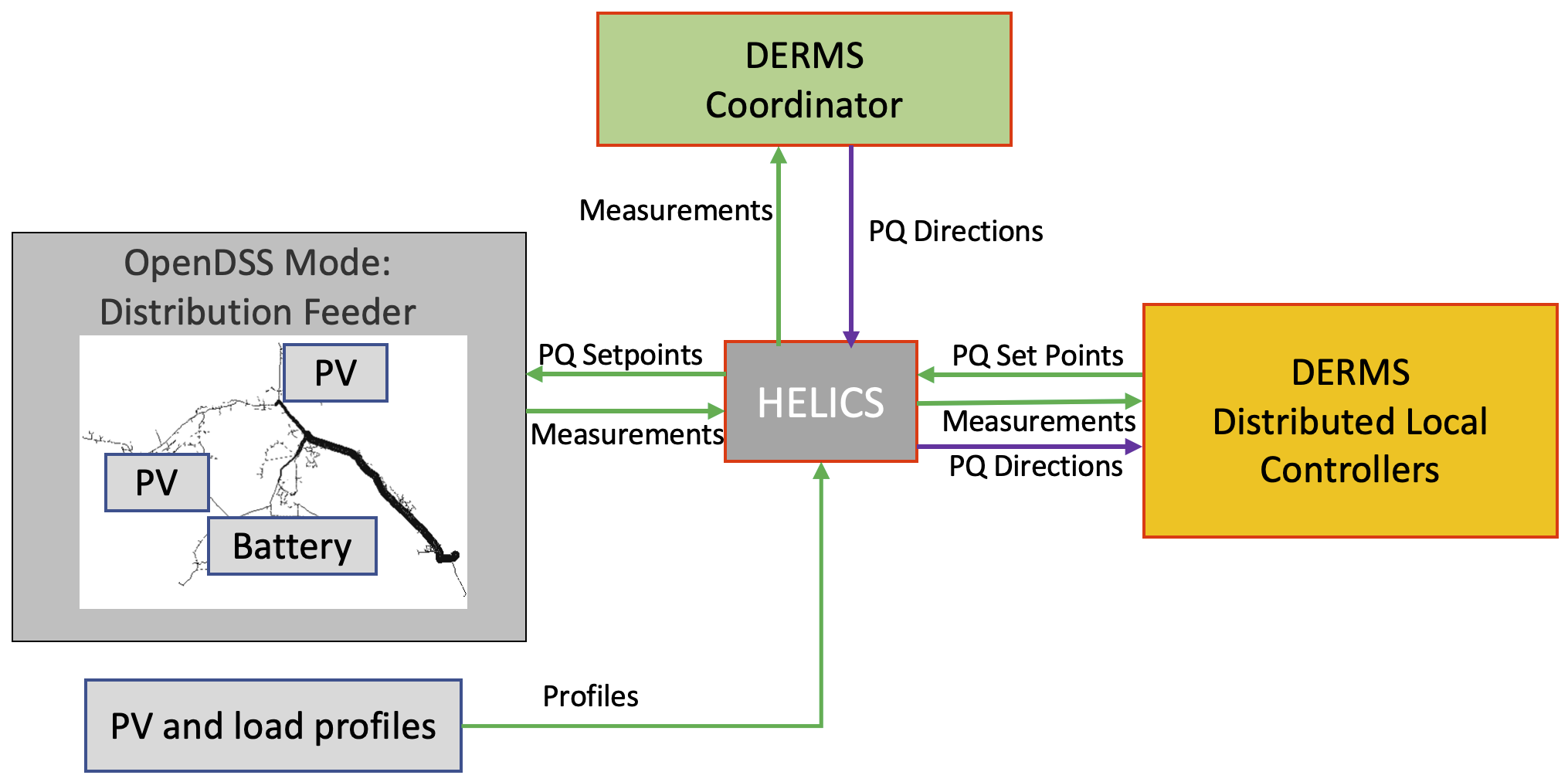}
    \caption{Co-simulation diagram with HELICS.}
    \label{fig:sim_setup_diagram}
\end{figure}

The system, which includes the distribution feeder, coordinator, and local controllers, is co-simulated using the Hierarchical Engine for Large-scale Infrastructure Co-Simulation (HELICS)~\cite{palmintier2017design}.
See Figure \ref{fig:sim_setup_diagram} for a diagram of the information being exchanged among the system components.
The real-world feeder is simulated as a quasi-steady-state time series in OpenDSS every 2 sec; it sends the voltage magnitude measurements and active power measurements at the feeder head to the coordinator and the implemented power injections of the DERs to their associated local controllers.
The coordinator sends the power injection direction signals to the local controllers every 2 sec, and the local controllers update their power injection set points every 2 sec.
The parameter settings for the primal-dual control algorithm are set to $\nu=10^{-3}$ and $\epsilon=10^{-4}$~\cite{dall2017optimal,bernstein2019real}.
The load and PV generation data were provided by a utility from their advanced metering infrastructure, with the loads changing every 15 min and the PV generation changing every 1 min.

The capacity of the PV generators range from 0.04 kW to 34 kW, with an average of 10 kW, and the batteries can store between 13.5 kWh and 54 kWh, with an average of 19 kWh.
Each of the 163 customers has a local controller that receives control signals from the DER management coordinator to decide the power set points of its PV generation and battery, if it has one.
The cost of curtailing the PV solar generator, $i$, from the available generation, $P_{\text{av},i,t}$, at time $t$ is $f_{i,t}(P_{i,t},Q_{i,t})=c^{(P)}_{i}(P_{i,t}-P_{\text{av},i,t})^2 + c^{(Q)}_{i}(Q_{i,t})^2$, where $c^{(P)}_{i}=\frac{0.2}{\text{INV}_i}$, $c^{(Q)}_{i}=\frac{0.002}{\text{INV}_i}$, and $\text{INV}_i$ is the rating of its inverter in kilowatts.
The feasible power injection set for the PV solar generator $i$ is $\mathcal{Y}_{i,t}=\{(P_{i,t},Q_{i,t}):0\leq P_{i,t}\leq P_{\text{av},i,t}, P_{i,t}^2+P_{i,t}^2\leq \text{INV}_i^2 \}$.
The cost of pushing the state of charge, $\text{SOC}_{i,t}$, of the battery, $i$, away from the preferred SOC, $\text{SOC}_{\text{pref},i}=60\%$, at time $t$ is $f_{i,t}(P_{i,t})=0.01(\frac{\text{SOC}_{i,t}}{100}-\frac{P_{i,t}dt}{\text{CAP}_i}-\frac{\text{SOC}_{\text{pref},i}}{100})^2$, where $\text{CAP}_i$ is the size of the battery in kWh, and $dt$ is the time between control decisions in hours.
The feasible active power injection for the battery, $i$, is $\mathcal{Y}_{i,t}=\{P_{i,t}:P_{\text{min},i,t}\leq P_{i,t} \leq P_{\text{max},i,t}\}$, where $P_{\text{min},i,t}$ and $P_{\text{max},i,t}$ are the minimum and maximum active power injections, respectively, which depend on the SOC and the charging/discharging limits of the battery.
These cost functions and feasible sets are described in more detail in \cite{padullaparti2022evaluation,bernstein2019real}.
The local controllers are also used as sensors measuring the voltage magnitudes and sending the values to the DER management coordinator.

\begin{figure}
    \centering
    \includegraphics[width=0.7\columnwidth]{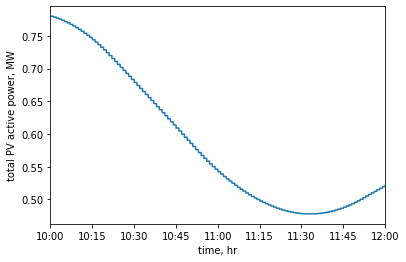}
    \caption{Total PV generation on 4/3/19 used in the default settings.}
    \label{fig:pv_total_active}
\end{figure}

The default parameter settings used in each testing scenario gives a baseline scenario that is stable and relatively static from which to make comparisons; they are the following, unless otherwise noted.
The initial step sizes are set to $\alpha_{i,0}=10:\forall i\in\mathcal{I}$, $\beta_{\text{volt},0}=5000$, and $\beta_{\text{VPP},0}=10$, which are relatively close to the values that they converge to when the default settings are used.
The initial step sizes are modified in Section \ref{sec:self-tuning}.
The VPP set points are set to the constant values of 1.02, 0.88, and 0.85 MW for phases A, B, and C, respectively, and they are changed to be time-varying in Section \ref{sec:vpp} and different constant values in Section \ref{sec:solar}.
The day and time are 4/3/19 from 10 am--12 pm, which has a smooth PV profile (see Figure \ref{fig:pv_total_active}), and they are changed in Section \ref{sec:solar} to have a more volatile PV profile.
The tap position of the LTC at the feeder head is set to -1 for all three phases, and it is changed to be time-varying in Section \ref{sec:tap_ratio}.

\subsection{Self-Tuning}
\label{sec:self-tuning}

During the process of manually tuning the step sizes, the operator will choose values that are either too low or too high before settling on a reasonable value.
With the adaptive procedure for the step sizes, the operator will need to choose only the initial step sizes for the adaptive procedure to start, and let the DER management tune itself from that point.
We showcase the DER management tuning itself from the initial step sizes $\alpha_{i,0}=0.1:\forall i\in\mathcal{I}$, $\beta_{\text{volt},0}=50$, and $\beta_{\text{VPP},0}=0.1$, which are two orders of magnitude less than the default settings and the values they converge to (labeled as ``low'').
This is repeated but with the initial step sizes $\alpha_{i,0}=1000:\forall i\in\mathcal{I}$, $\beta_{\text{volt},0}=5\times 10^5$, and $\beta_{\text{VPP},0}=1000$, which are two orders of magnitude greater than the default settings (labeled as ``high'').
Both scenarios are compared against the more optimistic default scenario, where the initial step sizes are close to what they converge to (labeled as ``base''). 

\begin{figure}
     \centering
     \begin{subfigure}[b]{0.6\columnwidth}
         \centering
         \includegraphics[width=\columnwidth]{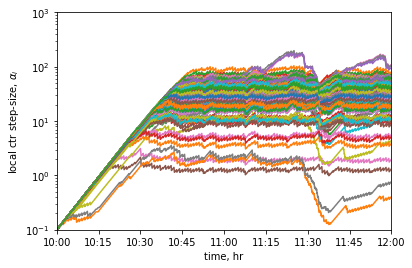}
         \caption{All step sizes start low}
         \label{fig:selftune_auto_low_stepsize_LCs}
     \end{subfigure}
     \hfill
     \begin{subfigure}[b]{0.6\columnwidth}
         \centering
         \includegraphics[width=\columnwidth]{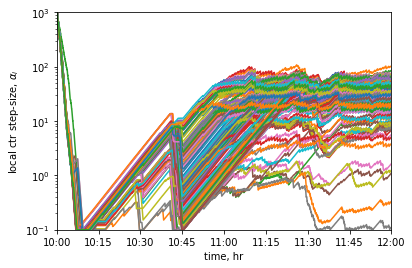}
         \caption{All step sizes start high}
         \label{fig:selftune_auto_high_stepsize_LCs}
     \end{subfigure}
     \hfill
     \begin{subfigure}[b]{0.6\columnwidth}
         \centering
         \includegraphics[width=\columnwidth]{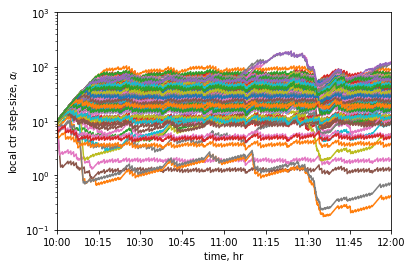}
         \caption{All step sizes start near converged values}
         \label{fig:selftune_auto_base_stepsize_LCs}
     \end{subfigure}
        \caption{Step size self-tuning of the local controllers.  Each color represents an individual controller, and they are consistent among (a), (b), and (c).}
        \label{fig:selftune_auto_stepsize_LCs}
\end{figure}

\begin{figure}
     \centering
     \begin{subfigure}[b]{0.7\columnwidth}
         \centering
         \includegraphics[width=\columnwidth]{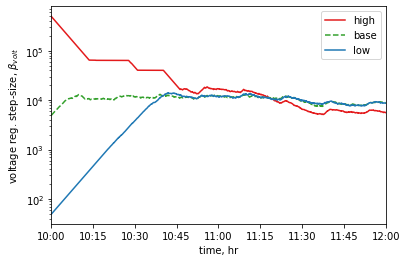}
         \caption{Voltage regulation}
         \label{fig:selftune_auto_stepsize_volt}
     \end{subfigure}
     \hfill
     \begin{subfigure}[b]{0.7\columnwidth}
         \centering
         \includegraphics[width=\columnwidth]{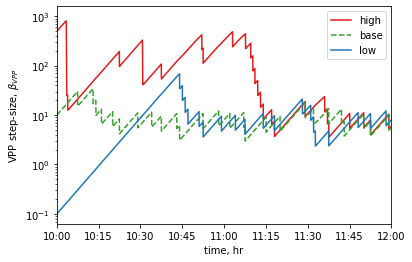}
         \caption{VPP}
         \label{fig:selftune_auto_stepsize_VPP}
     \end{subfigure}
        \caption{Step size self-tuning of the grid services.}
        \label{fig:selftune_auto_grid_services}
\end{figure}

The step sizes of the local controllers for the three scenarios are shown in Figure \ref{fig:selftune_auto_stepsize_LCs}, with each line color representing a specific local controller.
Although the initial step sizes vary by orders of magnitude between the scenarios, by the end of the two hours, they end up at very similar values, or at the very least similar magnitudes.
Figure \ref{fig:selftune_auto_grid_services} shows the step sizes of the grid services.
Like the step sizes of the local controllers, those of the voltage regulation and the VPP grid services also end up at very similar values.
This demonstrates the ability of the adaptive tuning procedure to reliably self-tune the DER management regardless of its initial conditions.

\begin{figure}
     \centering
     \begin{subfigure}[b]{0.6\columnwidth}
         \centering
         \includegraphics[width=\columnwidth]{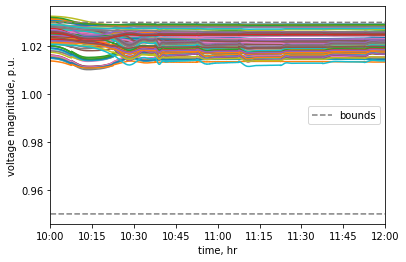}
         \caption{All step sizes start low.}
         \label{fig:selftune_auto_low_voltages}
     \end{subfigure}
     \hfill
     \begin{subfigure}[b]{0.6\columnwidth}
         \centering
         \includegraphics[width=\columnwidth]{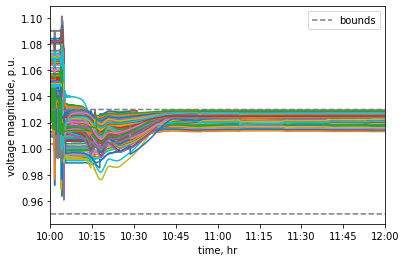}
         \caption{All step sizes start high.}
         \label{fig:selftune_auto_high_voltages}
     \end{subfigure}
     \hfill
     \begin{subfigure}[b]{0.6\columnwidth}
         \centering
         \includegraphics[width=\columnwidth]{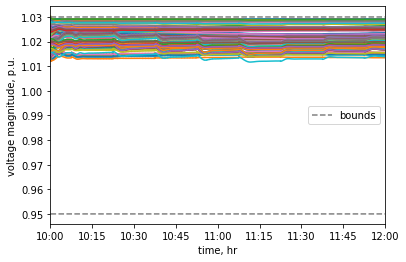}
         \caption{All step sizes start near converged values.}
         \label{fig:selftune_auto_base_voltages}
     \end{subfigure}
        \caption{Voltage magnitudes for the self-tuning scenarios.}
        \label{fig:selftune_auto_voltages}
\end{figure}

The most interesting aspect of the self-tuning property is that if the step sizes are the cause of instability in the system, i.e., they are too high and make the primal-dual DER control too reactive, the system can recover by quickly dropping the values of their step sizes.
For example, when the initial step sizes are extremely high, the voltages are far outside their bounds in the first 5 minutes compared to the scenario, when the initial step sizes are too low or near their converged values (see Figure \ref{fig:selftune_auto_voltages});
however, the voltage magnitudes are within bounds after 15 minutes and look very similar to the other scenarios within 45 minutes.
This shows that the self-tuning property can also be thought of as a self-correcting property, especially for the described extreme cases.

\subsection{VPP Set Point Tracking}
\label{sec:vpp}

The DER management operator providing the VPP grid service can be directed by a higher authority (e.g., ADMS) to change the set point of the feeder head active powers for a larger grid objective, e.g., frequency regulation for the transmission system.
Changes in the set point can be smooth or abrupt in a variety of ways; we evaluate step changes in the set point that are common and nonsmooth.
After the first 30 minutes, the set point of Phase A is increased by 0.26 MW, whereas the set points of phases B and C are decreased by 0.30 MW and 0.34 MW, respectively.
After another 45 minutes, the set point of Phase A is decreased by 0.13 MW, whereas the set points of phases B and C are increased by 0.15 MW and 0.17 MW, respectively.

\begin{figure}
     \centering
     \begin{subfigure}[b]{0.8\columnwidth}
         \centering
         \includegraphics[width=\columnwidth]{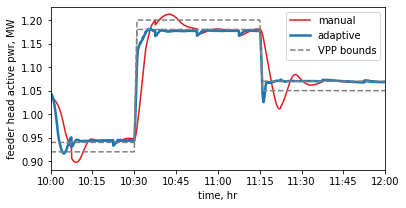}
         \caption{Phase A}
         \label{fig:vpptrack_PA}
     \end{subfigure}
     \hfill
     \begin{subfigure}[b]{0.8\columnwidth}
         \centering
         \includegraphics[width=\columnwidth]{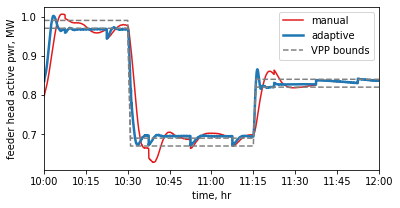}
         \caption{Phase B}
         \label{fig:vpptrack_PB}
     \end{subfigure}
     \hfill
     \begin{subfigure}[b]{0.8\columnwidth}
         \centering
         \includegraphics[width=\columnwidth]{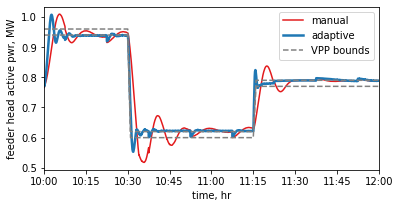}
         \caption{Phase C}
         \label{fig:vpptrack_PC}
     \end{subfigure}
        \caption{VPP set point tracking: manually tuned versus adaptively tuned step sizes.}
        \label{fig:vpptrack_P}
\end{figure}

\begin{figure}
     \centering
     \begin{subfigure}[b]{0.7\columnwidth}
         \centering
         \includegraphics[width=\columnwidth]{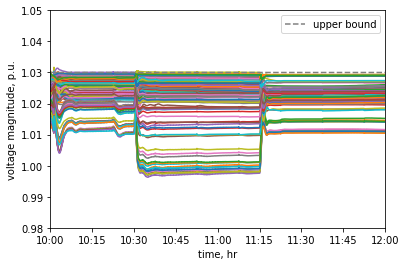}
         \caption{Adaptive step size tuning}
         \label{fig:vpptrack_auto_voltages}
     \end{subfigure}
     \hfill
     \begin{subfigure}[b]{0.7\columnwidth}
         \centering
         \includegraphics[width=\columnwidth]{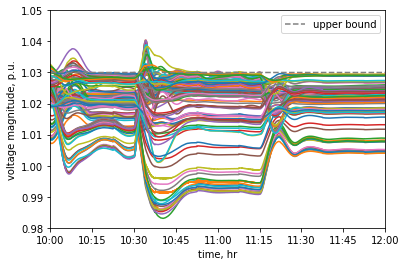}
         \caption{Manual step size tuning}
         \label{fig:vpptrack_man_voltages}
     \end{subfigure}
        \caption{Voltages magnitudes when the VPP set points have step changes.}
        \label{fig:vpptrack_voltages}
\end{figure}

\begin{figure}
     \centering
     \begin{subfigure}[b]{0.8\columnwidth}
         \centering
         \includegraphics[width=\columnwidth]{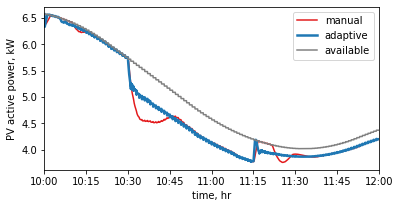}
         \caption{PV active power}
         \label{fig:vpptrack_localctr_PV_act}
     \end{subfigure}
     \hfill
     \begin{subfigure}[b]{0.8\columnwidth}
         \centering
         \includegraphics[width=\columnwidth]{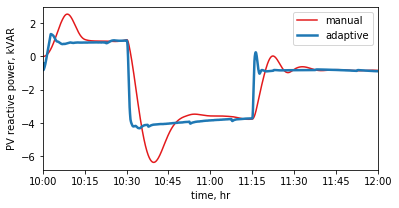}
         \caption{PV reactive power}
         \label{fig:vpptrack_localctr_PV_rea}
     \end{subfigure}
     \hfill
     \begin{subfigure}[b]{0.8\columnwidth}
         \centering
         \includegraphics[width=\columnwidth]{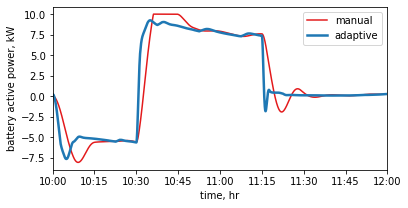}
         \caption{Battery active power}
         \label{fig:vpptrack_localctr_batt_act}
     \end{subfigure}
        \caption{Local controller example from Phase A: manually tuned versus adaptively tuned step sizes.}
        \label{fig:vpptrack_localctr}
\end{figure}

The VPP set point tracking for all three phases is shown in Figure \ref{fig:vpptrack_P} and compares the use of adaptively tuned versus manually tuned step sizes.
At each step change, the adaptively tuned step sizes have a faster and less oscillatory transition than the manually tuned ones.
The measured voltage magnitudes are shown in Figure \ref{fig:vpptrack_voltages}, where the adaptively tuned step sizes keep the voltage within bounds with fewer oscillations than the manually tuned step sizes.
Also, the voltage magnitudes of the manually tuned step sizes spike 0.01 p.u. above the upper bound, whereas those of the adaptively tuned step sizes do not.
The cleaner transitions of the adaptively tuned step sizes are also observed at the local controller level.
Figure \ref{fig:vpptrack_localctr} gives an example of the power injections of the PV generation and battery for a local controller that confirms this observation.

\begin{figure}
     \centering
     \begin{subfigure}[b]{0.7\columnwidth}
         \centering
         \includegraphics[width=\columnwidth]{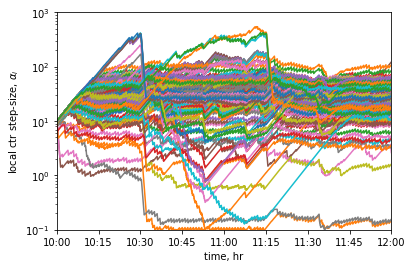}
         \caption{Local controllers}
         \label{fig:vpptrack_auto_stepsize_LCs}
     \end{subfigure}
     \hfill
     \begin{subfigure}[b]{0.7\columnwidth}
         \centering
         \includegraphics[width=\columnwidth]{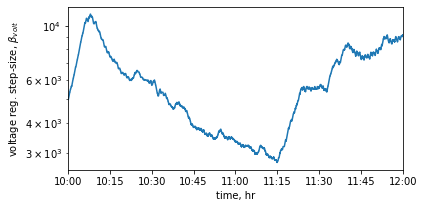}
         \caption{Voltage regulation}
         \label{fig:vpptrack_auto_stepsize_volt}
     \end{subfigure}
     \begin{subfigure}[b]{0.7\columnwidth}
         \centering
         \includegraphics[width=\columnwidth]{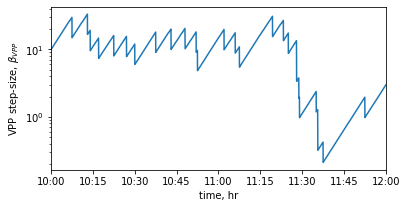}
         \caption{VPP}
         \label{fig:vpptrack_auto_stepsize_VPP}
     \end{subfigure}
        \caption{Step sizes of the local controllers and the grid services during step changes in the VPP set points.}
        \label{fig:vpptrack_auto_stepsize}
\end{figure}

The step sizes for the local controller and grid services are shown in Figure \ref{fig:vpptrack_auto_stepsize}.
The step change in the VPP set point at 10:30 causes no observable change on its step size and instead causes changes in the step sizes of the local controllers.
The step change at 11:15, however, causes both the VPP step size and the local controllers to change.
The determination of which step sizes change depends on where in the system the oscillations are observed via the independent use of the cosine similarity at the system components.

\subsection{Fluctuating PV Generation}
\label{sec:solar}

\begin{figure}
    \centering
    \includegraphics[width=0.7\columnwidth]{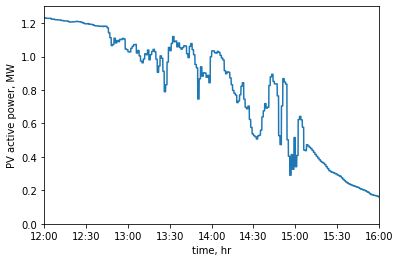}
    \caption{Total available PV solar generation.}
    \label{fig:solar_total_pv}
\end{figure}

One main concern when adding large amounts of PV generation capacity to a distribution network is that there are days when the amount of available power fluctuates greatly, which can significantly stress the control system. 
For example, a sudden unexpected drop in PV generation caused by a cloud passing by will cause both the voltage magnitude to drop and the power at the feeder head to increase abruptly; a DER management operator will need to quickly counteract this.
We test this with a PV profile that starts near its peak generation at noon and then has large sudden fluctuations between 12:30 and 15:00 (see Figure \ref{fig:solar_total_pv}).
For the testing of this scenario, we increase the priority of the VPP by increasing the step size tuning parameter, $\underline{\gamma}_\text{VPP}$, to 0.995 and by decreasing the DER parameters, $\underline{\gamma}_i:\forall i\in\mathcal{I}$, to 0.8 from the default settings.

\begin{figure}
     \centering
     \begin{subfigure}[b]{0.8\columnwidth}
         \centering
         \includegraphics[width=\columnwidth]{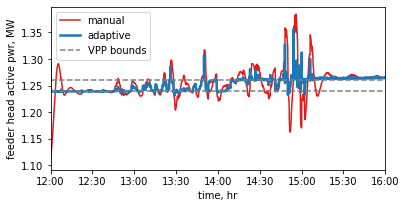}
         \caption{Phase A}
         \label{fig:solar_PA}
     \end{subfigure}
     \hfill
     \begin{subfigure}[b]{0.8\columnwidth}
         \centering
         \includegraphics[width=\columnwidth]{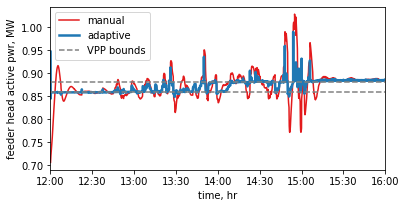}
         \caption{Phase B}
         \label{fig:solar_PB}
     \end{subfigure}
     \hfill
     \begin{subfigure}[b]{0.8\columnwidth}
         \centering
         \includegraphics[width=\columnwidth]{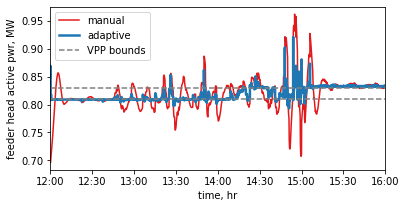}
         \caption{Phase C}
         \label{fig:solar_PC}
     \end{subfigure}
        \caption{VPP under fluctuating PV solar generation: manually tuned versus adaptively tuned step sizes.}
        \label{fig:solar_P}
\end{figure}

\begin{figure}
     \centering
     \begin{subfigure}[b]{0.7\columnwidth}
         \centering
         \includegraphics[width=\columnwidth]{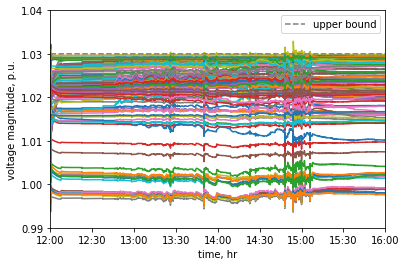}
         \caption{Adaptive step size tuning}
         \label{fig:solar_auto_voltages}
     \end{subfigure}
     \hfill
     \begin{subfigure}[b]{0.7\columnwidth}
         \centering
         \includegraphics[width=\columnwidth]{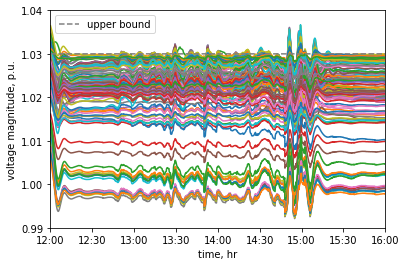}
         \caption{Manual step size tuning}
         \label{fig:solar_man_voltages}
     \end{subfigure}
        \caption{Voltages magnitudes when the available PV solar generation fluctuates significantly.}
        \label{fig:solar_voltages}
\end{figure}

The VPP set point trajectories are shown in Figure \ref{fig:solar_P}.
The adaptively tuned step sizes are better able to keep the feeder head powers within the VPP bounds and prevent them from oscillating as much as the manually tuned ones;
however, it is not possible to prevent them from spiking above the VPP upper bound.
This is caused by the available PV generation suddenly dropping below its set point and instantaneously cutting generation power from the feeder head.
This requires at least a full measurement and decision iteration of the DERMS to be counteracted.
When the available PV generation suddenly increases, its local controller can gradually increase the set point and avoid suddenly dropping the feeder head power below the VPP bound.
The voltage magnitudes in Figure \ref{fig:solar_voltages} show that the adaptively tuned step sizes prevent oscillations and better respect the bounds compared to the manually tuned ones.

\begin{figure}
     \centering
     \begin{subfigure}[b]{0.7\columnwidth}
         \centering
         \includegraphics[width=\columnwidth]{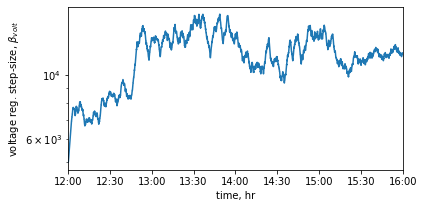}
         \caption{Voltage regulation}
         \label{fig:solar_auto_stepsize_volt}
     \end{subfigure}
     \begin{subfigure}[b]{0.7\columnwidth}
         \centering
         \includegraphics[width=\columnwidth]{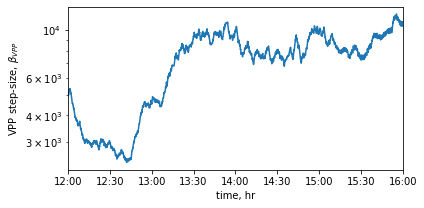}
         \caption{VPP}
         \label{fig:solar_auto_stepsize_VPP}
     \end{subfigure}
        \caption{Step sizes of the local controllers and the grid services during step changes in the VPP set points.}
        \label{fig:solar_auto_stepsize}
\end{figure}

The adaptively tuned step sizes are given in Figure \ref{fig:solar_auto_stepsize} for the voltage regulation and VPP grid services.
Even though the PV generation fluctuates, it causes the step sizes for the grid services to increase instead of decreasing because the fluctuations add more pressure on the control system to keep the services within their bounds.

\subsection{Tap Ratio Changes}
\label{sec:tap_ratio}

\begin{figure}
    \centering
    \includegraphics[width=0.7\columnwidth]{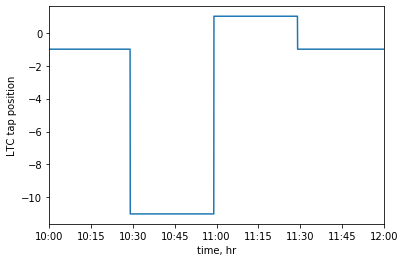}
    \caption{Tap position of the LTC.}
    \label{fig:ltc_tap_position}
\end{figure}

When a higher-level controller (e.g., ADMS) makes a decision that abruptly changes the admittance matrix of a distribution network, such as changing the tap ratio, it can have two major effects on the DERMS.
The first is that it might instantaneously make the currently implemented DER power injections not satisfy one or more of the grid services.
This will cause the coordinator to send control signals that are significantly different than the previous one and abruptly change the DER power injections to satisfy the grid services.
In this case, if the step sizes associated with the grid services are too large, there can be large oscillations from reacting to the sudden change.
The second major effect is that there can be a mismatch between the admittance matrix of the distribution network and the one assumed by the DER management operator if the change was not communicated to the operator.
This can cause the control signals sent by the coordinator to be suboptimal with respect to satisfying the grid services.

In this numerical simulation, we assume that an ADMS is also tasked with keeping the voltage magnitudes between bounds by changing the tap position of the LTC located at the feeder head.
Because of the shared responsibility of the voltage regulation grid service between the ADMS and DER management operator, we lower the priority of the voltage regulation on the DER management side by decreasing the parameter $\underline{\gamma}_\text{volt}$ to 0.25 from its default setting.
We implement large step changes in both directions of the tap position, which are shown in Figure \ref{fig:ltc_tap_position}, to simulate the ADMS trying to manipulate the voltage magnitudes to be within bounds, but it overshoots by pushing them too low at 10:30 and then pushing them too high at 11:00.

\begin{figure}
     \centering
     \begin{subfigure}[b]{0.7\columnwidth}
         \centering
         \includegraphics[width=\columnwidth]{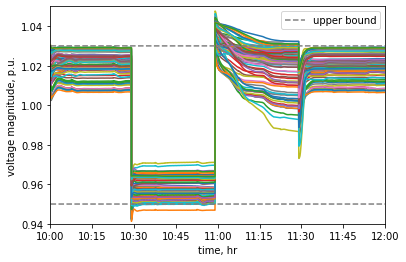}
         \caption{Adaptive step size tuning}
         \label{fig:ltc_auto_voltages}
     \end{subfigure}
     \hfill
     \begin{subfigure}[b]{0.7\columnwidth}
         \centering
         \includegraphics[width=\columnwidth]{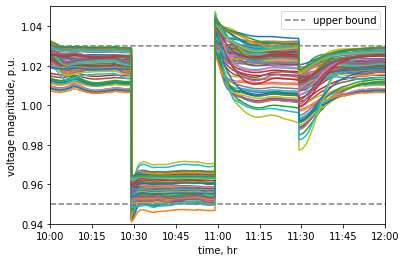}
         \caption{Manual step size tuning}
         \label{fig:ltc_man_voltages}
     \end{subfigure}
        \caption{Voltages magnitudes when the LTC tap position changes.}
        \label{fig:ltc_voltages}
\end{figure}

\begin{figure}
     \centering
     \begin{subfigure}[b]{0.8\columnwidth}
         \centering
         \includegraphics[width=\columnwidth]{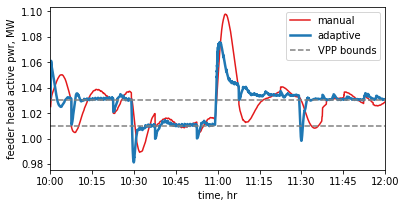}
         \caption{Phase A}
         \label{fig:ltc_PA}
     \end{subfigure}
     \hfill
     \begin{subfigure}[b]{0.8\columnwidth}
         \centering
         \includegraphics[width=\columnwidth]{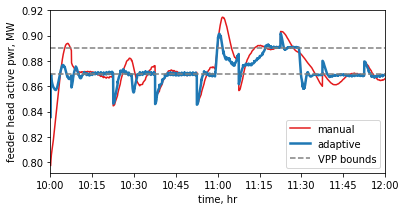}
         \caption{Phase B}
         \label{fig:ltc_PB}
     \end{subfigure}
     \hfill
     \begin{subfigure}[b]{0.8\columnwidth}
         \centering
         \includegraphics[width=\columnwidth]{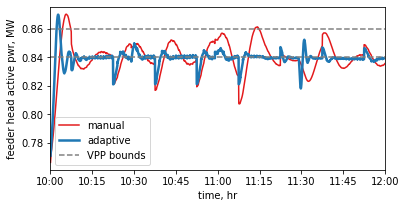}
         \caption{Phase C}
         \label{fig:ltc_PC}
     \end{subfigure}
        \caption{VPP under changes in the LTC tap position: manually tuned versus adaptively tuned step sizes.}
        \label{fig:ltc_P}
\end{figure}

The voltage magnitudes are displayed in Figure \ref{fig:ltc_voltages}, which shows them instantaneously spiking down when the tap position is dropped and spiking up when the tap position is raised.
In both cases, primal-dual DER control gradually brings them back closer to their bounds;
however, the immediate reaction to the voltage magnitudes suddenly being outside of their bounds has an adverse impact on the VPP grid service.
Specifically, in phases A and B, the feeder head active powers (see Figure \ref{fig:ltc_P}) spike below their bounds at 10:30 and spike above them at 11:00, which are caused by the DER management primarily reacting to the voltage magnitudes being shifted outside of their bounds; the stronger the immediate reaction is to the voltage bound violations, the larger the spike is in the feeder head powers.
Afterward, the shifts in the feeder head active power being pushed outside of their bounds cause the DER control to react to both grid services simultaneously in almost opposite directions with respect to the DER power injections.
In fact, the gap above the feeder head active power bound for Phase A and the gap above the voltage magnitude bound from 11:15--11:30 is caused by the two grid services pushing against each other, neither being completely satisfied.

The final observation is that the feeder head powers with adaptively tuned step sizes have fewer and smaller oscillations than the manually tuned ones because of their ability to counteract oscillations wherever they are measured in the system.

\subsection{Discussion}

\begin{table}[]
    \centering
    \begin{tabular}{|l|l|} 
    \hline
    Scenario & Notes  \\
    \hline\hline
    Self-tuning & The step sizes either gradually increase to make the \\
    & DER management more responsive when they are \\
    & too low or quickly decrease to avoid oscillations when \\
    & they are too high.  \\ 
    \hline
    VPP set point & DER management can track the VPP set point more  \\
    Tracking & precisely and with fewer oscillations. \\
    \hline
    Fluctuating PV & DER management has tighter voltage regulation and \\
    generation & better prevents VPP oscillations below its lower \\
    & bound. \\
    \hline
    Tap ratio & DER management reacts more aggressively to keep \\
    changes & the VPP within bounds. \\
    \hline
    \end{tabular}
    \caption{Performance summary of using adaptive step size tuning over manual step size tuning}
    \label{tab:perf}
\end{table}

We have demonstrated that adaptive step size tuning has many advantages over manually tuning the step sizes beyond the main motivation that it removes a labor-intensive step for an engineer.
Table \ref{tab:perf} gives a summary of its advantages based on the evaluated scenarios.

The main general advantage, as demonstrated by the last three scenarios, is that it can perform more aggressively and precisely than its manually tuned version.
This is because the manually tuned version only allows for a global step size to be tuned, which must be small enough not to cause oscillations in any part of the system, thus limiting how aggressive it can be.
The adaptive tuning framework, however, splits the step sizes to be individualized to each control component so that a component requiring a low step size does not limit the size of another.
The adaptive tuning quickly removes any oscillations that appear at any component so that each component can be as aggressive as possible without causing instability.

Finally, the self-tuning scenario shows that the performance of adaptive tuning does not depend on its initial settings when given a reasonable amount of time, which means that it can leave one control scenario that it adapted to and adapt to another.

%% file: sections/conclusion.tex
In this paper, we designed a simple method for automatically and adaptively tuning the step sizes of a DERMS that uses primal-dual control.
Based on the conditions of the distribution network, the method will either increase the step sizes to accelerate the reaction of the DERs to provide grid services more quickly, or it will decrease the step sizes to remove any observed oscillations.
This adaptive tuning procedure reduces the amount of time that an engineer is required to spend setting up the DERMS for deployment.
Additionally, the method allows the DER management operator to prioritize among possibly competing grid services.
The method was evaluated using numerical simulations of a real-world feeder in Colorado using data obtained from an electric utility.
We showed that the primal-dual DERMS using the automatically tuned step sizes can better adapt to changing VPP bounds, highly fluctuating PV solar power, and changes in LTC tap positions than a DERMS with a step size that was manually tuned.